\documentclass{amsart}
\usepackage{amsfonts}
\usepackage{amsmath,amssymb}
\usepackage{amsthm}
\usepackage{amscd}
\usepackage{graphics}
\usepackage{graphicx}
\theoremstyle{definition}{
	\newtheorem{Def}{{\rm Definition}}
	\newtheorem{Ex}{{\rm Example}}
	
	\newtheorem{Prob}{{\rm Problem}}
}
\theoremstyle{plain}
{
	
	\newtheorem{Prop}{Proposition}
	\newtheorem{Thm}{Theorem}
	\newtheorem{MainThm}{Main Theorem}

}

\begin{document}
	\title[A class of naturally generalized special generic maps]{A class of naturally generalized special generic maps}
	\author{Naoki Kitazawa}
	\address{Institute of Mathematics for Industry, Kyushu University, 744 Motooka, Nishi-ku Fukuoka 819-0395, Japan\\
		TEL (Office): +81-92-802-4402 \\
		FAX (Office): +81-92-802-4405 \\
	}
	\email{n-kitazawa@imi.kyushu-u.ac.jp}
	\urladdr{https://naokikitazawa.github.io/NaokiKitazawa.html}
	
	\begin{abstract}
		{\it Special generic} maps are generalizations of Morse functions with exactly two singular points on spheres and canonical projections of unit spheres. They restrict the manifolds of the domains strongly in considerable cases and are important in algebraic topology and differential topology of manifolds of specific classes and manifolds regarded as elementary in some senses admit such maps in considerable cases.
		
		We propose a class of generalized special generic maps in our present paper and extend a fundamental result on structures and some algebraic topological properties of special generic maps by the author. Our present study will be a pioneering study on nice classes of generalized special generic maps. Studies of algebraic topological properties and differential topological ones of special generic maps have developed due to their nice structures for example.


		
	\end{abstract}
	
	
	\maketitle
	\section{Introduction.}
	\label{sec:1}
	
Morse functions are important smooth functions, discussed systematically in \cite{milnor1,milnor2} for example. Their {\it singular} points, or points where the ranks of the differentials drop, appearing discretely, know information on the homology groups and some homotopy of the manifolds.
Morse functions with exactly two singular points characterizing spheres topologically except $4$-dimensional cases. This is Reeb's theorem.

Special generic maps generalize the functions and canonical projections of spheres embedded in the Euclidean spaces in canonical ways, or {\it unit spheres}.
\\
By its definition, they restrict the topologies and the differentiable structure of the manifolds strongly, whereas several manifolds regarded as elementary in some senses. They can cover wide classes of manifolds (of the domains) when we restrict the manifolds as ones seen as elementary in some senses admit such maps in considerable cases. Pioneering studies are \cite{burletderham, furuyaporto} for example, \cite{calabi} studies related problems as another pioneering one. \cite{saeki1, saeki2, saekisakuma1, saekisakuma2} are pioneering studies related to algebraic topology and differential topology of manifolds. They are followed by \cite{nishioka,wrazidlo1,wrazidlo2} for example. They mainly study special generic maps on spheres, manifolds represented as connected sums of the products of (two) spheres and restrictions on the homology groups of the manifolds.
Related to these studies,  the author has started studies on the cohomology rings and he has shown that the real projective spaces and complex projective spaces admit no such maps unless they are $1$-dimensional real or complex projective spaces, in other words, circles or $2$-dimensional spheres for example: see \cite{kitazawa1, kitazawa4, kitazawa5} and for related studies see also \cite{kitazawa3, kitazawa6} for example. {\it Non-singular} {\it compact complex toric} varieties or {\it toric} manifolds are important in complex algebraic geometry and (equivariant) differential topology. They seem to admit no such maps in considerable cases. Complex projective spaces are of such a class. For toric manifolds, see \cite{buchstaberpanov} for example.

We can say that special generic maps can cover manifolds of some specific classes as domains well and that special generic maps cannot cover considerable large classes of manifolds. This presents the following problem.

\begin{Prob}
	\label{prob:1}
	What are good classes of smooth maps whose codimensions are not positive which cover manifolds being the domains of some special generic maps and manifolds not being so such as general toric manifolds?
 
\end{Prob}
Special generic maps have very simple structures from the viewpoint of algebraic topology and differential topology. Due to this, studies on algebraic topological studies and differential topological ones of the maps and the manifolds have developed. Another problem is as follows.
\begin{Prob}
	\label{prob:2}
	For classes of Problem \ref{prob:1}, do maps there have such nice properties?
\end{Prob} 
 
 In this paper, we do not give direct answers to them. Related to this, we introduce such smooth maps as {\it simply generalized special generic maps} generalizing the class of special generic maps and study some fundamental or explicit properties which will be fundamental, meaningful and strong in studies in the future.
 More precisely, the image of a special generic map on a closed and connected manifolds is a smoothly immersed compact and connected manifold whose codimension is $0$. The preimage of a point in the interior of the immersed manifold is regarded as a sphere (before it is immersed). Around the boundary it is represented as the product of a natural Morse function ({\it height function}) on a disk and the identity map on a disk and in the interior it is regarded as a projection. Here, Morse functions are replaced by so-called {\it Morse-Bott} functions and the preimage of a point in the interior is a product of spheres. For Morse-Bott functions see \cite{bott} for example.
 
We present one of our Main Theorems.
We leave some undefined terminologies, notions and notation later. Here we introduce several terminologies, notions and notation on differentiable manifolds and maps. Elementary notions from algebraic topology such as (co)homology groups, cohomology rings, and elements of homology groups represented by compact and connected oriented (sub)manifolds with no boundaries are presented later or consult \cite{hatcher} for example. 

The $k$-dimensional Euclidean space is a simplest smooth manifold and it is also a Riemannian manifold equipped with the so-called standard Euclidean metric. Let it be denoted by ${\mathbb{R}}^k$ for $k \geq 1$ and $\mathbb{R}:={\mathbb{R}}^1$. $\mathbb{Z} \subset \mathbb{R}$ denotes the set of all integers, regarded as a commutative ring which is a principal ideal domain having the identity element $1$ different from the zero element $0$ naturally. For each $x \in {\mathbb{R}}^k$ $||x|| \geq 0$ denotes the distance between this and the origin $0 \in {\mathbb{R}}^k$ under the metric. $S^k:=\{x \in {\mathbb{R}}^{k+1} \mid ||x||=1 \}$ denotes the {\it $k$-dimensional unit sphere} for $k \geq 0$, which is a $k$-dimensional smooth compact submanifold of ${\mathbb{R}}^{k+1}$ with no boundary. This is connected for $k \geq 1$ and the two-point discrete set for $k=0$. $D^k:=\{x \in {\mathbb{R}}^{k} \mid ||x|| \leq 1 \}$ denotes the {\it $k$-dimensional unit disk} for $k \geq 1$, which is a $k$-dimensional smooth compact and connected submanifold of ${\mathbb{R}}^{k+1}$ and whose boundary is the unit sphere $S^{k-1}$.

\begin{MainThm}[Some of Theorem \ref{thm:3}]
\label{mthm:1}
Let $n$ be an arbitrary positive integer. Let ${\bar{f}}_N:\bar{N} \rightarrow N$ be a smooth immersion of an $n$-dimensional compact, connected and orientable manifold $\bar{N}$ into an $n$-dimensional connected and orientable manifold $N$ with no boundary.
Let $l_1>0$ be an integer and let $m_{l_1}$ be a map from the set of all positive integers smaller than or equal to $l_1$ to the set of all positive integers.
Let $\mathcal{C}_{\bar{N}}:=\{C_{j-1}\}_{j=1}^{l_2}$ denote the set of all connected components of the boundary $\partial \bar{N} \subset \bar{N}$, consisting of exactly $l_2 \geq 0$ connected components. 
Let $m_{\mathcal{C}_{\bar{N}},l_1}$ be a map from the set before into the set of all positive integers smaller than or equal to $l_1$.

Then we have a simply generalized special generic map $f_0:M_0 \rightarrow N$ on a suitable closed and connected manifold $M_0$ of dimension $m:=n+{\Sigma}_{j=1}^{l_1} m_{l_1} (j)$ enjoying the following properties.
	\begin{enumerate}
		\item $f_0:M_0 \rightarrow N$ is represented as the composition of a smooth surjection $q_{f_0}:M_0 \rightarrow \bar{N}$ with the given immersion ${\bar{f}}_N$. The preimage ${q_{f_0}}^{-1}(p)$ is diffeomorphic to the product ${\prod}_{j=1}^{l_1} S^{m_{l_1}(j)}$ for any $p \in {\rm Int}\ \bar{N}$.
 		\item For any principal ideal domain $A$ having a unique identity element different from the zero element and $H_i(M_0;A)$ is free for any integer $i$, the following properties are enjoyed.
\begin{enumerate} 
		\item The homology group $H_i(M_0;A)$ has a submodule isomorphic to $H_i(\bar{N};A)$ seen as a module over $A$ for any integer $i$.
		\item The cohomology ring $H^{\ast}(M_0;A)$ has a subalgebra isomorphic to $H^{\ast}(\bar{N};A)$ where they are seen as graded commutative algebra over $A$.
		\item For each integer $1 \leq j \leq l_2-1$, We have the submodule of the module $H_{i-1}(M_0;A)$ generated by the finite set of elements each of which is represented by a submanifold in $M_0$ diffeomorphic to the product of $C_j$ and the product ${\prod}_{j^{\prime}=1}^{l_1} {S^{\prime}}_{P,B,C_j,f_0,j^{\prime}}$, defined in the following way for any integer $i>n$.
		
		\begin{itemize}
			\item We have the product ${\prod}_{j^{\prime}=1}^{l_1} S_{P,B,C_j,f_0,j^{\prime}}$ of standard spheres or one-point sets satisfying the following conditions.
			\begin{itemize}
				\item $S_{P,B,C_j,f_0,j^{\prime}}:=S^{m_{l_1}(j^{\prime})}$ if $j^{\prime} \notin m_{{\mathcal{C}}_{\bar{N}},l_1}(\{C_0,C_j\})$.
			\item $S_{P,B,C_j,f_0,j^{\prime}}$ is a one-point set in $S^{m_{l_1}(j^{\prime})}$ if $j^{\prime} \in m_{{\mathcal{C}}_{\bar{N}},l_1}(\{C_0,C_j\})$.
			\end{itemize}
			\item For the product ${\prod}_{j^{\prime}=1}^{l_1} {S^{\prime}}_{P,B,C_j,f_0,j^{\prime}}$,
            ${S^{\prime}}_{P,B,C_j,f_0,j^{\prime}}$ is ${S}_{P,B,C_j,f_0,j^{\prime}}$ or a one-point set there. Furthermore, this product is not a one-point set.
			
		\end{itemize}
		
\end{enumerate}.
		\end{enumerate}
\end{MainThm}

This is related to Problem \ref{prob:3}, presented later. This gives an answer to an question asking what is a simplest simply generalized special generic map whose image is a given smoothly immersed manifold and what is the topology of the manifold of the domain. The case of special generic maps has been solved by the author in the preprint \cite{kitazawa8} and also presented in ones \cite{kitazawa2, kitazawa7}. Theorem \ref{thm:1} presents this again. 

What follows is the organization of our paper. The second section is for preliminaries. We present several fundamental terminologies, notions and notation on complexes including polyhedra, differentiable maps, and diffeomorphisms and bundles. After that we review special generic maps. Generalized versions of special generic maps are introduced as our new work and main objects and tools in our paper. We also review fundamental classes of compact, connected (oriented) manifolds and elements of homology groups represented by such (sub)manifolds. We define our new class rigorously in the third section. The fourth section is on Main Theorems. 

	\ \\
	{\bf Conflict of interest.} \\
	The author is a member of the project JSPS KAKENHI Grant Number JP22K18267 "Visualizing twists in data through monodromy" (Principal Investigator: Osamu Saeki). Our present study is supported by the project. \\
	\ \\
	{\bf Data availability.} \\
	Data essentially supporting our present study are all contained in our present paper.
	
	\section{Preliminaries.}
	We review important terminologies, notions, notation and properties.
\subsection{A short exposition on cell complexes, CW complexes, polyhedra and PL manifolds.}
For a topological space $X$ having the structure of a cell complex the maximal dimension for cells among cells of which is finite, we can define the dimension uniquely, denoted by $\dim X$.

(Topological) manifolds are known to be given the structures of CW complexes. Every smooth manifold is canonically regarded as a polyhedron, which is a so-called {\it PL manifold}. For each smooth manifold, this canonically obtained PL manifold is defined uniquely.

According to so-called Hauptvermutung, discussed in \cite{moise} for example, a topological manifold whose dimension is at most $3$ has a PL structure uniquely and it does not have PL structures other than this.

For example, \cite{hudson} is important for general theory on the PL category and the piecewise smooth category, which are known to be equivalent, especially on PL topological theory. We omit rigorous expositions on fundamental or advanced theory here.
\subsection{Differentiable maps and singular points.}
Let $c:X \rightarrow Y$ be a differentiable map between from a differentiable manifold $X$ into $Y$. $p \in X$ is a {\it singular} point of $c$ if the rank of the differential ${dc}_p$ is smaller than both $\dim X$ and $\dim Y$ at $p$ and $c(p)$ is a {\it singular value} of $c$. Let $S(c)$ denote the set of all singular points of $c$; this is the {\it singular set} of $c$.
\subsection{Diffeomorphisms and bundles.}
	A {\it diffeomorphism} between smooth manifolds is a smooth map with no singular points which is also a homeomorphism. A {\it diffeomorphism on a manifold} means a diffeomorphism from the (smooth) manifold onto itself. Two manifolds are said to be {\it diffeomorphic} if and only if there exists a diffeomorphism between the two manifolds. This naturally gives an equivalence relation on the family of all smooth manifolds where their corners are eliminated in a canonical way. 
	Note also that this canonical way gives a smooth manifold with no corner and if we consider these operations to a fixed manifold, then we always have manifolds which are mutually diffeomorphic.
	We can define {\it PL homeomorphic manifolds} using piecewise smooth homeomorphisms similarly.
	
	The {\it diffeomorphism group} of a manifold means the space consisting of all diffeomorphisms on the manifold endowed with the so-called {\it Whitney $C^{\infty}$ topology} and it is of course a group. 
	Note that the {\it Whitney $C^{\infty}$ topology} on the space of all smooth maps between two smooth manifolds and its subspaces are a natural and important topology and important spaces in the theory of singularity theory of differentiable maps and differential topology of manifolds. See \cite{golubitskyguillemin}.
		
	A {\it smooth} bundle means a bundle whose fiber is a smooth manifold and whose structure group is the diffeomorphism group.
	A {\it linear} bundle means a bundle whose fiber is a Euclidean space, a unit disk, or a unit sphere, and whose structure group consists of linear transformations where linear transformations can be defined naturally and diffeomorphisms of course.
	
	Precise explanation on general theory of bundles is omitted and see \cite{steenrod} for example. For linear bundles, see \cite{milnorstasheff} for example. 
	
	\subsection{Special generic maps.}
	\begin{Def}
		\label{def:1}
		A smooth map $c:X \rightarrow Y$ between two smooth manifolds with no boundaries are said to be a {\it special generic} map if at each singular point $p \in X$ we have suitable local coordinates around $p$ (and $f(p)$) and for them we can represent $c$ as $(x_1,\cdots,x_{\dim X}) \rightarrow (x_1,\cdots,x_{\dim Y-1},{\Sigma}_{j=1}^{\dim X-\dim Y+1} {x_{\dim Y+j-1}}^2)$.
	\end{Def}
	 A canonical projection of a unit sphere is defined as a map mapping $(x_1,x_2) \in S^k \subset {\mathbb{R}}^{k+1}={\mathbb{R}}^{k_1} \times {\mathbb{R}}^{k_2}$ to $x_1 \in {\mathbb{R}}^{k_1}$ with $k \geq 1$, $k_1,k_2 \geq 1$ and $k=k_1+k_2$. To show that this is a special generic map is an exercise on smooth maps, Morse functions and singularity theory of differentiable maps. 
	 
	 Since the 1990s, manifolds admitting special generic maps have been actively studied by Saeki and Sakuma. Related results are in \cite{saeki1,saeki2,saekisakuma1,saekisakuma2}, Following these studies, Nishioka and Wrazidlo have obtained  \cite{nishioka,wrazidlo1,wrazidlo2}.
	They have revealed restrictions on the differentiable structures of the homotopy spheres, some elementary manifolds such as ones in Example \ref{ex:1}, which is presented later. Restrictions on the homology groups have been also studied. As a pioneer, the author launched explicit systematic studies on the cohomology rings of the manifolds in \cite{kitazawa1,kitazawa2,kitazawa3,kitazawa4,kitazawa5,kitazawa6,kitazawa7,kitazawa8} for example.

	\begin{Prop}[\cite{saeki1,saeki2}]
		\label{prop:1}
		Given a special generic map $f:M \rightarrow N$ on an $m$-dimensional closed and connected manifold $M$ into an $n$-dimensional connected manifold $N$ which has no boundary with $m \geq n \geq 1$. This enjoys the following properties.
		\begin{enumerate}
			\item \label{prop:1.1}
			We have a suitable $n$-dimensional compact and connected smooth manifold $W_f$, a suitable smooth surjection $q_f:M \rightarrow W_f$ and a suitable smooth immersion $\bar{f}:W_f \rightarrow N$ in such a way that they enjoy the relation $f=\bar{f} \circ q_f$. Furthermore, we can have $q_f$ as a map whose restriction to the singular set $S(f)$ of $f$ is a diffeomorphism onto the boundary $\partial W_f \subset W_f$.
			\item \label{prop:1.2}
Let the singular set $S(f)$ be non-empty. Then we can have a small collar neighborhood $N(\partial W_f)$ of the boundary $\partial W_f \subset W_f$ in such a way that the following two properties are enjoyed.
			\begin{enumerate}
				\item \label{prop:1.2.1}
				The composition of the restriction of $q_f$ to the preimage ${q_f}^{-1}(N(\partial W_f))$ with the canonical projection to $\partial W_f$ is the projection of a linear bundle whose fiber is the unit disk $D^{m-n+1}$.
				\item \label{prop:1.2.2}
				The restriction of $q_f$ to the preimage of $W_f-{\rm Int}\ N(\partial W_f)$ is the projection of a smooth bundle whose fiber is the unit sphere $S^{m-n}$. Furthermore, in some specific case, we can do so that the bundle is linear and the case $m-n=0,1,2,3$ gives such an example.
		\end{enumerate}
	\end{enumerate}
	\end{Prop}

\begin{Def}[E. g. \cite{kitazawa8}]
	\label{def:2}
	In Proposition \ref{prop:1}, we call the bundle of (\ref{prop:1.2.1}) the {\it boundary linear bundle} of $f$. There we call the bundle of (\ref{prop:1.2.2}) the {\it internal smooth bundle} of $f$.
\end{Def}

The following gives simplest special generic maps.

	\begin{Prop}[\cite{saeki1}]
		\label{prop:2}
Let $m \geq n \geq 1$ be integers. Let $\bar{N}$ be an $n$-dimensional smooth, compact and connected manifold whose boundary $\partial \bar{N}$ is not empty and $N$ an $n$-dimensional smooth connected manifold which has no boundary. Assume also that a smooth immersion ${\bar{f}}_N:\bar{N} \rightarrow N$ is given.
			
		Then we have a suitable $m$-dimensional closed and connected manifold $M$ and some special generic map $f:M \rightarrow N$ enjoying the following properties.
		\begin{enumerate}
			\item The property {\rm (}\ref{prop:1.1}{\rm )} of Proposition \ref{prop:1} where $W_f$ and $\bar{N}$ are identified in a suitable way
                        \item The property that the relation ${\bar{f}}_N=\bar{f}$ is enjoyed.
			\item A boundary linear bundle and an internal smooth bundle of $f$ are both trivial.
		\end{enumerate}
 
	\end{Prop}

We present a main result of \cite{kitazawa8}, which is also presented in \cite{kitazawa2, kitazawa7} for example in suitably revised forms.
For notions and theorems from elementary algebraic topology such as the {\it cup product} for an ordered pair or a sequence of a finite length of elements of cohomology groups, the intersection form for a compact and connected (oriented) manifold, derived naturally from so-called Poincar\'e duality (theorem) for the manifold, and universal coefficient theorem, see \cite{hatcher} again for example,
\begin{Thm}[\cite{kitazawa8}]
	\label{thm:1}
	A special generic map in Proposition \ref{prop:2} can be constructed as a map $f_0:M_0 \rightarrow N$ enjoying the following properties on a suitably chosen manifold $M_0:=M$ where $f_0:=f$ and $W_{f_0}:=W_f$ for example in the proposition.
	\begin{enumerate}
		\item \label{thm:1.1}
 We have a special generic map or a smooth immersion $f_j:M_0 \rightarrow N \times {\mathbb{R}}^{j}$ for any positive integer $j$ such that for any pair of distinct non-negative integers $j_1<j_2$, $f_{j_1}={\pi}_{j_2,j_1,N} \circ f_{j_2}$ holds for the canonical projection ${\pi}_{j_2,j_1,N}:N \times {\mathbb{R}}^{n+j_2} \rightarrow N \times {\mathbb{R}}^{n+j_1}$ with $0<j_1<j_2$ or ${\pi}_{j_2,0,N}:N \times {\mathbb{R}}^{n+j_2} \rightarrow N$ with $j_2>j_1=0$. 
		\item \label{thm:1.2}
 For any principal ideal domain $A$ having the unique identity element different from the zero element, the homology group such that $H_j(W_{f_0};A)$ is free as a module over $A$ for any integer $j$, the following properties are enjoyed.
\begin{enumerate} 
		\item The homology group $H_i(M_0;A)$ is isomorphic to the direct sum of $H_i(W_{f_0};A)$ and $H_{i-(m-n)}(W_{f_0},\partial W_{f_0};A)$ as a module over $A$ and we identify the direct summand suitably as submodules of $H_i(M_0;A)$ and do so for the cohomology groups for any integer $i${\rm :} note that universal coefficient theorem yields similar representations.
		\item The cohomology ring $H^{\ast}(W_{f_0};A)$ gives a subalgebra of the cohomology ring $H^{\ast}(M_0;A)$ isomorphic to this subalgebra where they are seen as graded commutative algebras over $A$.
		
		\item The cup product $u_1 \cup u_2$ of the ordered pair of an element $u_1 \in H^{i_1}(W_{f_0};A) \subset H^{j_1}(M_0;A)$ and an element $u_2 \in H^{i_2-(m-n)}(W_{f_0};A) \subset H^{i_2}(M_0;A)$ gives a bilinear form over $A$ isomorphic to the intersection form for the manifold $W_{f_0}$ if it is oriented.
		\item 
		The cup product of elements in $\oplus H^{i-(m-n)}(W_{f_0},\partial W_{f_0};A) \subset \oplus H^i(M_0;A)$ is always the zero element except for the case $i=0$ with $m=n$ where "$\oplus$" is for the direct sums for every integer $i$.
		\end{enumerate}
\end{enumerate}
\end{Thm}

	Other than canonical projections of unit spheres, we present simplest special generic maps.
	\begin{Ex}
		\label{ex:1}
		Let $l$ be an arbitrary positive integer and $m \geq n \geq 2$ integers. We give an integer $1 \leq n_j \leq n-1$ for each integer $1 \leq j \leq l$. We take a connected sum of $l>0$ manifolds considered in the smooth category where the $j$-th manifold is $S^{n_j} \times S^{m-n_j}$ and have a smooth manifold $M_0$. We have a special generic map $f_0:M_0 \rightarrow {\mathbb{R}}^n$ as in Proposition \ref{prop:2} and Theorem \ref{thm:1}. More precisely, we have $f_0$ in such a way that the image is represented as a boundary connected sum of $l$ manifolds considered in the smooth category with the $j$-th manifold being diffeomorphic to $S^{n_j} \times D^{n-n_j}$ 
		
	\end{Ex}
	
    Hereafter, a {\it homotopy sphere} means a smooth manifold homeomorphic to a unit sphere whose dimension is not $0$. 
	A {\it standard} sphere means a homotopy sphere being diffeomorphic to some unit sphere where an {\it exotic sphere} means a homotopy sphere which is not. It is well-known that $4$-dimensional exotic spheres are still unknown. Except these cases, homotopy spheres are known to be PL homeomorphic to unit spheres as canonically defined PL manifolds. $4$-dimensional exotic spheres are known to be not PL homeomorphic to standard spheres where they are considered to be the canonical PL manifolds.
	
	\begin{Thm}[\cite{saeki1,saeki2}]
		\label{thm:2}
		\begin{enumerate}
			\item
			\label{thm:2.1}
			Let $m$ be an arbitrary integer greater than or equal to $2$.
			A necessary and sufficient condition for an $m$-dimensional closed and connected manifold $M$ to admit a special generic map $f:M \rightarrow {\mathbb{R}}^2$ is that $M$ is either of the following manifolds.
			\begin{enumerate}
				\item A homotopy sphere which is not a $4$-dimensional exotic sphere.
				\item A manifold
				represented as a connected sum of smooth manifolds taken in the smooth category where each manifold is the total space of a smooth bundle over $S^1$ whose fiber is either of the following two.
				\begin{enumerate}
					\item An {\rm (}$m-1${\rm )}-dimensional homotopy sphere with $m \neq 5$.
					\item A $4$-dimensional standard sphere.
				\end{enumerate}
			\end{enumerate}
			\item
			\label{thm:2.2}
			Let $m$ be an arbitrary integer greater than or equal to $4$.
				We have a sufficient condition for an $m$-dimensional closed and simply-connected manifold $M$ admitting a special generic map $f:M \rightarrow {\mathbb{R}}^3$. If $M$ is either of the following manifolds, then $M$ admits one.
			\begin{enumerate}
			\item A homotopy sphere which is not a $4$-dimensional exotic sphere.
				\item A manifold
				represented as a connected sum of smooth manifolds taken in the smooth category where each manifold here is the total space of a smooth bundle over $S^2$ whose fiber is either of the following two.
				\begin{enumerate}
				\item An {\rm (}$m-2${\rm )}-dimensional homotopy sphere where $m \neq 6$.
	      		\item A $4$-dimensional standard sphere.
				\end{enumerate}
			\end{enumerate}
			In the case $m=4,5,6$, the converse is also true. In such a case, a fiber of each bundle is an {\rm (}$m-2${\rm)}-dimensional standard sphere and the total spaces of the bundles are regarded the total spaces of linear bundles without changing the fibers and the base spaces.
			\item
			\label{thm:2.3}
		Both in the cases {\rm (}\ref{thm:2.1}{\rm )} and {\rm (}\ref{thm:2.2}{\rm )},
		the manifold $M$ not being a homotopy sphere admits a special generic map as in Example \ref{ex:1} where the internal smooth bundle and the boundary linear bundle of it are not trivial in general. In the two cases, let $n_j=1$ and $n_j=2$ in Example \ref{ex:1}, respectively.
		\end{enumerate}
	\end{Thm}

\subsection{Fundamental classes of connected and compact  (oriented) manifolds.}

For terminologies, notions and notation on (co)homology groups and systematic explanation related to arguments here, see \cite{hatcher} again for example.

Let $A$ be a commutative ring having a unique identity element different from the zero element and let $1_A$ denote the identity element. $1_A$ and $-1_A$ are generators of $A$ where $A$ is considered to be the canonically obtained module over $A$.
For a compact and connected oriented manifold $X$, $H_{\dim X}(X, \partial X; A)$ is isomorphic to $A$ as a module over $A$. We can give a generator according to the orientation of $X$. In related arguments, we do not need orientations in the case $A:=\mathbb{Z}/2\mathbb{Z}$ or the commutative ring of order $2$.

For a manifold $Y$, consider an embedding $i_X:X \rightarrow Y$ with suitable conditions being assumed according to the category where we discuss. For example, in the smooth category, this is smooth and in the PL or the piecewise smooth category, this is piecewise smooth. Furthermore, $X$ is embedded {\it properly}, in other words, the boundary is embedded into the boundary and the interior is embedded into the interior and (in the smooth category) $X$ must be embedded in a so-called generic way.  
If $a \in H_j(Y, \partial Y; A)$ is the value of the homomorphism ${i_{X}}_{\ast}:H_{\dim X}(X, \partial X; A) \rightarrow H_{\dim X}(Y, \partial Y; A)$ canonically induced from the embedding at the fundamental class $[X] \in H_{\dim X}(X, \partial X; A)$, then $a$ is said to be {\it represented} by the submanifold $i_X(X)$.

We note on a "generic way" in embedding manifolds smoothly. In short it satisfies conditions on so-called "transversality". 
This is important in several scenes of our paper.
For our smooth embedding $i_X:X \rightarrow Y$ here, we consider an embedding satisfying the condition that 
the dimension of the intersection of the image of the differential $d {i_X}_p$ of the embedding $i_X:X \rightarrow Y$ at each point $p \in \partial X$ in the boundary and the tangent space at $i_X(p) \in \partial Y$ is $\dim X+\dim \partial Y-\dim Y=\dim X-1$. This is also important in singularity theory of differentiable maps and differential topology. For such notions and arguments, see \cite{golubitskyguillemin} again.

\section{Our new class of smooth maps--simply generalized special generic maps--.}

A {\it height function of a unit sphere} is a function which is a specific case of canonical projections of unit spheres. This is a Morse function of course.
This is defined as a map mapping $x=(x_1,x_2) \in S^k \subset {\mathbb{R}}^{k+1}=\mathbb{R} \times {\mathbb{R}}^k$ to $x_1$.
A {\it height function of a unit disk} is defined as the restriction of the height function of the unit sphere to the preimage of $\mathbb{R} \bigcap \{t \geq 0 \mid t \in \mathbb{R}\}$.

\begin{Def}
	\label{def:3}
	Let $m \geq n \geq 1$ be integers.
	A {\it simply generalized special generic} map is a smooth map $f:M \rightarrow N$ on an $m$-dimensional closed and connected manifold $M$
	 such that the image is a smooth immersion ${\bar{f}}_N:\bar{N} \rightarrow N$ of an $n$-dimensional compact and connected manifold $\bar{N}$ into an $n$-dimensional smooth manifold $N$ which has no boundary and that enjoys the following properties.
	 \begin{enumerate}
	 	\item \label{def:3.1}
 As in Proposition \ref{prop:1}, we have a smooth surjection $q_f:M \rightarrow W_f$ enjoying the relation $f={\bar{f}}_N \circ q_f$ where $\bar{N}$ is identified in a suitable way with $W_f$ as a smooth manifold. Furthermore, we can define $q_f$ as a map whose restriction to the singular set defines a diffeomorphsm onto the boundary $\partial W_f$.
	 	\item \label{def:3.2}
 We have some small collar neighborhood $N(\partial W_f)$ for the boundary $\partial W_f\subset W_f$ and the composition of the restriction of $q_f$ to the preimage with the canonical projection to the boundary is the projection of a smooth bundle over $\partial W_f$.
\item \label{def:3}
 On the collar neighborhood $N(\partial W_f)$ and the preimage ${q_f}^{-1}(N(\partial W_f))$, it is locally represented as a smooth function represented as the product map of the following two smooth maps for local suitable coordinates around each point of $\partial W_f \subset N(\partial W_f)$. 
\begin{enumerate}
\item \label{def:3.3.1}
 The composition of the projection of a trivial smooth bundle over a unit disk whose fiber is either of the following two and a height function on the unit disk.
\begin{enumerate}
\item
\label{def:3.3.1.1}
 The product of finitely many standard spheres.
\item
\label{def:3.3.1.2}
 The disjoint union of copies of a manifold as in (\ref{def:3.3.1.1}) of an even number.
\end{enumerate}
Furthermore, the domain of this function is regarded as a fiber of the smooth bundle over $\partial W_f$ just before.
\item
\label{def:3.3.2}
 The identity map on a small smoothly embedded copy of the unit disk $D^{\dim \partial W_f} \subset \partial W_f$.
\end{enumerate}
 
	 	\item
\label{def:3.4}
 The restriction of $q_f$ to the preimage ${q_f}^{-1}(W_f-{\rm Int}\ N(\partial W_f))$ is the projection of a smooth bundle whose fiber is
either of the following two.

\begin{enumerate}
\item
\label{def:3.4.1.1}
 The product of finitely many standard spheres.
\item
\label{def:3.4.1.2}
 The disjoint union of copies of a manifold as in (\ref{def:3.4.1.1}) of an even number .
\end{enumerate}

	 \end{enumerate}
	
\end{Def}

We present a fundamental proposition without its proof or we can easily check this. 

\begin{Prop}
Special generic maps are seen as specific simply generalized special generic maps if their singular sets are not empty.
\end{Prop}

In short, it is a smooth map locally represented as the projection or the product map of a so-called {\it Morse-Bott} function of a good class and the identity map on a disk. More precisely, Morse-Bott functions here are for each component of a so-called {\it moment map} on a {\it symplectic toric} manifold of dimension $2k$ into ${\mathbb{R}}^k$ for example. For related theory, see \cite{buchstaberpanov, delzant} for example.

\section{On Main Theorems.}

We concentrate on the following natural and fundamental problem.

\begin{Prob}
\label{prob:3}
	For a given smooth immersion ${\bar{f}}_N:\bar{N} \rightarrow N$ whose codimension is $0$ as in Proposition \ref{prop:2}, Theorem \ref{thm:1} and Definition \ref{def:3} for example, can we construct a simplest generalized special generic map whose image is as the given smoothly immersed manifold. Theorem \ref{thm:1} gives an answer for special generic maps.
	In addition, can we investigate algebraic topological properties and differential topological ones of the map and the resulting manifold of the domain as Theorem \ref{thm:1} shows.
\end{Prob}

The following gives an answer to this and should be regarded as key facts in challenging Problems \ref{prob:1} and \ref{prob:2}. 

\begin{Thm}[Some are in Main Theorem \ref{mthm:1}]
\label{thm:3}
Let $n$ be an arbitrary positive integer. Let ${\bar{f}}_N:\bar{N} \rightarrow N$ be a smooth immersion of an $n$-dimensional compact, connected and orientable manifold $\bar{N}$ into an $n$-dimensional connected and orientable manifold $N$ with no boundary.
Let $l_1>0$ be an integer and let $m_{l_1}$ be a map from the set of all positive integers smaller than or equal to $l_1$ to the set of all positive integers.
Let $\mathcal{C}_{\bar{N}}:=\{C_{j-1}\}_{j=1}^{l_2}$ denote the set of all connected components of the boundary $\partial \bar{N} \subset \bar{N}$, consisting of exactly $l_2 \geq 0$ connected components. 
Let $m_{\mathcal{C}_{\bar{N}},l_1}$ be a map from the set before into the set of all positive integers smaller than or equal to $l_1$.

Then we have a simply generalized special generic map $f_0:M_0 \rightarrow N$ on a suitable closed and connected manifold $M_0$ of dimension $m:=n+{\Sigma}_{j=1}^{l_1} m_{l_1} (j)$ enjoying the following properties.
\begin{enumerate}
	\item \label{thm:3.1}
	$f_0:M_0 \rightarrow N$ is represented as the composition of a smooth surjection $q_{f_0}:M_0 \rightarrow W_{f_0}$ with the given immersion ${\bar{f}}_N$ where $W_{f_0}$ and $\bar{N}$ are suitably identified as smooth manifolds.
	\item \label{thm:3.2}
	 The preimage ${q_{f_0}}^{-1}(p)$ is diffeomorphic to the product ${\prod}_{j=1}^{l_1} S^{m_{l_1}(j)}$ for $p \in {\rm Int}\ W_{f_0}$.
	\item \label{thm:3.3}
	For any principal ideal domain $A$ having the unique identity element $1_A$ different from the zero element $0_A$ and $H_i(M_0;A)$ is free for any integer $i$, the following properties are enjoyed.
	\begin{enumerate} 
		\item \label{thm:3.3.1}
		 The homology group $H_i(M_0;A)$ has a submodule isomorphic to $H_i(W_{f_0};A)$ seen as a module over $A$ for any integer $i$. Let the submodule be denoted by ${H_i(W_{f_0};A)}_{H_i(M_0;A)} \subset H_i(M_0;A)$.
		\item \label{thm:3.3.2}
		For any integer $i$, we can choose some basis of ${H_i(W_{f_0};A)}_{H_i(M_0;A)} \subset H_i(M_0;A)$ and for any element of it, we have an element of $H_{m-i}(M_0;A)$ such that the cup product for the Poincar\'e duals to these two are a generator of $H^{m}(M_0;A)$ {\rm (}where the two are ordered suitably{\rm )}. Furthermore, such elements form a basis of a submodule of $H_{m-i}(M_0;A)$. Let the submodule be denoted by ${H_{n-i}(W_{f_0},\partial W_{f_0};A)}_{H_{m-i}(M_0;A)} \subset H_{m-i}(M_0;A)$.
		\item \label{thm:3.3.3}
		The intersection of ${H_i(W_{f_0};A)}_{H_i(M_0;A)} \subset H_i(M_0;A)$ and ${H_{n-m+i}(W_{f_0},\partial W_{f_0};A)}_{H_{i}(M_0;A)} \subset H_{i}(M_0;A)$ is the trivial module for any integer $i$.
				\item  \label{thm:3.3.4}
				 The cohomology ring $H^{\ast}(M_0;A)$ has a subalgebra isomorphic to $H^{\ast}(W_{f_0};A)$ where they are seen as graded commutative algebra over $A$.
			\item  \label{thm:3.3.5}
			For each integer $1 \leq j \leq l_2-1$, there exists a submodule of the module $H_{i-1}(M_0;A)$ generated by the set of elements each of which is represented by some {\rm (}$i-1${\rm )}-dimensional smooth submanifold in $M_0$ diffeomorphic to the product of $C_j$ and the product of ${\prod}_{j^{\prime}=1}^{l_1} {S^{\prime}}_{P,B,C_j,f_0,j^{\prime},a}$, defined in the following way for any integer $i$ satisfying the conditions $i>n$ and $i<m$.
		
		\begin{itemize}
			\item We have the product ${\prod}_{j^{\prime}=1}^{l_1} S_{P,B,C_j,f_0,j^{\prime}}$ of standard spheres or one-point sets satisfying the following conditions.
			\begin{itemize}
				\item $S_{P,B,C_j,f_0,j^{\prime}}:=S^{m_{l_1}(j^{\prime})}$ if $j^{\prime} \notin m_{{\mathcal{C}}_{\bar{N}},l_1}(\{C_0,C_j\})$.
				\item $S_{P,B,C_j,f_0,j^{\prime}}$ is a one-point set in $S^{m_{l_1}(j^{\prime})}$ if $j^{\prime} \in m_{{\mathcal{C}}_{\bar{N}},l_1}(\{C_0,C_j\})$.
			\end{itemize}
		
			\item 
			For the product ${\prod}_{j^{\prime}=1}^{l_1} {S^{\prime}}_{P,B,C_j,f_0,j^{\prime},a}$,
			${S^{\prime}}_{P,B,C_j,f_0,j^{\prime},a}$ is ${S}_{P,B,C_j,f_0,j^{\prime}}$ or a one-point set there. Furthermore, this product is not a one-point set.
			
		\end{itemize}
	Moreover, the product of $C_j$ and the product of ${\prod}_{j^{\prime}=1}^{l_1} {S}_{P,B,C_j,f_0,j^{\prime}}$ can be regarded as a smooth submanifold with no boundary in $M_0$ in a suitable natural way. Furthermore, the product of $C_j$ and the product of ${\prod}_{j^{\prime}=1}^{l_1} {S^{\prime}}_{P,B,C_j,f_0,j^{\prime},a}$ is also regarded as its submanifold in a canonical way. 
	
	Under this situation, an element represented by the product of $C_{j}$ and the product of ${\prod}_{j^{\prime}=1}^{l_1} {S^{\prime}}_{P,B,C_{j},f_0,j^{\prime},a}$ is not the zero element.
	Furthermore, an element represented by the product of $C_{j_1}$ and the product of ${\prod}_{j^{\prime}=1}^{l_1} {S^{\prime}}_{P,B,C_{j_1},f_0,j^{\prime},a_1}$ and an element represented by the product of $C_{j_2}$ and the product of ${\prod}_{j^{\prime}=1}^{l_1} {S^{\prime}}_{P,B,C_{j_2},f_0,j^{\prime},a_2}$ are mutually linearly independent if and only either of the following holds.
	\begin{itemize}
		\item $j_1 \neq j_2$.
		\item $j_1=j_2$ holds and the pair $({\prod}_{j^{\prime}=1}^{l_1} {S^{\prime}}_{P,B,C_{j_1},f_0,j^{\prime},a_1},{\prod}_{j^{\prime}=1}^{l_1} {S^{\prime}}_{P,B,C_{j_1},f_0,j^{\prime},a_2})$ is a pair of distinct subspaces of the same spaces ${\prod}_{j^{\prime}=1}^{l_1} S_{P,B,C_{j_1},f_0,j^{\prime}}={\prod}_{j^{\prime}=1}^{l_1} S_{P,B,C_{j_2},f_0,j^{\prime}}$.
	\end{itemize}
	\item \label{thm:3.3.6}
	For any submodule ${H_{n-i_1}(W_{f_0},\partial W_{f_0};A)}_{H_{m-i_1}(M_0;A)} \subset H_{m-i_1}(M_0;A)$ defined in {\rm (}\ref{thm:3.3.2}{\rm )} and any submodule of $H_{i_2-1}(M_0;A)$ with $i_2>n$ and $m-i_1=i_2-1$ presented in {\rm (}\ref{thm:3.3.5}{\rm )}, their intersection is always the trivial submodule.
	\end{enumerate}.
\end{enumerate}
\end{Thm}
We prove Theorem \ref{thm:3}. Theorem \ref{thm:1} can be regarded as a result for a specific case except (\ref{thm:1.1}) there. We review some of our proof of Theorem \ref{thm:1} in \cite{kitazawa8} in our present proof.  

Several theorems on elementary algebraic topology such as the (co)homology exact sequence for a pair of topological spaces and Mayer-Vietoris sequences, universal coefficient theorem and Künneth formula, which we need in the proof, are also in \cite{hatcher} for example. Let $A$ be a principal ideal domain $A$ having the unique identity element $1_A$ different from the zero element $0_A$.
We define the {\it cohomology dual} ${e_j}^{\ast,{\mathcal{E}}_{X,A}} \in H^i(X;A)$ to an element in a basis ${\mathcal{E}}_{X,A}=\{e_j\}$ respecting the basis of the unique maximal free submodule of a finite rank of the module $H_i(X;A)$. This can be uniquely defined as an element satisfying the relation ${e_{j_1}}^{\ast, {\mathcal{E}}_{X,A}}(e_{j_2})={\delta}_{j_1,j_2}$ where ${\delta}_{j_1,j_2}$ is $1_A$ if $j_1=j_2$ and $0_A$ if $j_1 \neq j_2$.

For a compact, connected and oriented manifold $X$, the {\it Poincar\'e dual} to an element of $H_j(X;A)$ is uniquely defined as an element of $H^{\dim X-j}(X,\partial X;A)$. Similarly, the {\it Poincar\'e dual} to an element of $H^j(X;A)$ is uniquely defined as an element of $H_{\dim X-j}(X,\partial X;A)$. Furthermore, the {\it Poincar\'e dual} to an element of $H_j(X,\partial X;A)$ is uniquely defined as an element of $H^{\dim X-j}(X;A)$ and the {\it Poincar\'e dual} to an element of $H^j(X,\partial X;A)$ is uniquely defined as an element of $H_{\dim X-j}(X;A)$.

Let $t_0,t_1 \in \mathbb{R}$. Let ${\delta}_{\geq t_0,t_1}:\mathbb{R} \rightarrow \mathbb{R}$ denote the function such that ${\delta}_{\geq t_0,t_1}(x)=0$ for $x<t_0$ and that ${\delta}_{\geq t_0,t_1}(x)=t_1$ for $x \geq t_0$. We also use such a function.

\begin{proof}[A proof of Theorem \ref{thm:3}]

We can construct a simply generalized special generic map $f_0:M_0 \rightarrow N$ on a suitable closed and connected manifold $M_0$ in the following way where our new smooth manifold  $W_{f_0}$ is identified with $\bar{N}$ as a smooth manifold in a suitable way. In addition, in Definition \ref{def:3} "$f$" is replaced by $f_0$ as $q_{f_0}$ and $W_{f_0}$ show for example. 
We construct a trivial smooth bundle over $W_{f_0}-{\rm Int}\ N(\partial W_{f_0})$ whose fiber is diffeomorphic to ${\prod}_{j_1=1}^{l_1} S^{m_{l_1}(j_1)}$. Let $P_{f_0}$ denote the total space.

$C_j \subset \partial W_{f_0}$ is a connected component of $\partial W_{f_0}$ and we define $N(C_j)$ as a connected component of $N(\partial W_{f_0})$ and a small collar neighborhood of $C_j$.
For each $C_j$, we glue the product map of the following maps along $(W_f-{\rm Int}\ N(\partial W_{f_0})) \bigcap \partial N(C_j)$ by the product map of the identity map between the bases spaces and the product map of the identity maps on the standard spheres where we need exposition on the identifications. Here the base spaces are identified suitably and naturally. In addition, standard spheres of the products of the standard spheres of the fibers can be given identifications suitably and naturally.
\begin{itemize}
\item The composition of the projection of a trivial smooth bundle over $D^{m_{l_1} \circ m_{{\mathcal{C}}_{\bar{N}},l_1}(C_j)+1}$ whose fiber is the product ${\prod}_{j_1=1}^{l_1-1} S^{m_{l_1}(j_1+{\delta}_{\geq m_{{\mathcal{C}}_{\bar{N}},l_1}(C_j),1}(j_1))}$ with a height function on $D^{m_{l_1} \circ m_{{\mathcal{C}}_{\bar{N}},l_1}(C_j)+1}$.
\item The identity map on $C_j$.
\end{itemize}

Let the domain of each map denoted by $B_{C_j,f_0}$ and the disjoint union by $B_{C,f_0}$. 

This finishes the presentation on the simply generalized special generic map $f_0$. This completes the proof of the properties (\ref{thm:3.1}) and (\ref{thm:3.2}).

We show the property (\ref{thm:3.3}). The assumption implies that homology groups and submodules of them here are free as modules over $A$ unless otherwise stated. 
Hereafter, let $l_2>0$ where we can show similarly in the case $l_2=0$. 

We show (\ref{thm:3.3.1}), (\ref{thm:3.3.2}), (\ref{thm:3.3.3}), and (\ref{thm:3.3.4}).

We can define a basis $\{e_{i,j}\}_{j \in J}$ of $H_i(W_{f_0};A)$ and for any element $e_{i,j} \in H_i(W_{f_0};A)$ here and the cohomology dual ${e_{i,j}}^{\ast, {\mathcal{E}}_{X,A}} \in H^i(W_{f_0};A)$ to it respecting the basis, there exists the Poincar\'e dual ${\rm PD}({e_{i,j}}^{\ast, {\mathcal{E}}_{X,A}}) \in H_{n-i}(W_{f_0},\partial W_{f_0};A)$ to it. 
We consider the fibers of the trivial smooth bundle over $W_{f_0}-{\rm Int}\ N(\partial W_{f_0})$ whose fiber is the product of finitely many standard spheres. The fiber of each trivial bundle over $(W_{f_0}-{\rm Int}\ N(\partial W_{f_0})) \bigcap \partial N(C_j)$ obtained by the restriction of the bundle is the boundary of the product of finitely many standard spheres and a copy of some unit disk and this bundle is a subbundle of a trivial smooth bundle whose fiber is the product of finitely many standard spheres and the copy of the unit disk. From this, we have an element ${\rm PD}({e_{i,j}}^{\ast, {\mathcal{E}}_{X,A}}) \in H_{n-i+{\sum}_{j=1}^{l_1} m_{l_1}(j)}(M_0;A)=H_{m-i}(M_0;A)$ where the same notation is used. In short, as presented in the proof of Theorem \ref{thm:1} in \cite{kitazawa8}, we consider variants of {\it Thom classes} of linear bundles or {\it prism operators}.

We have a section of the trivial bundle over $W_{f_0}-{\rm Int}\ N(\partial W_{f_0})$. This implies that the cohomology ring $H^{\ast}(M_0;A)$ has a subalgebra isomorphic to $H^{\ast}(W_{f_0};A)$ and $H^{\ast}(W_{f_0}-{\rm Int}\ N(\partial W_{f_0});A)$ where $W_{f_0}$ collapses to $W_{f_0}-{\rm Int}\ N(\partial W_{f_0})$, diffeomorphic to the original manifold $W_{f_0}$.

In addition, $H_i(M_0;A)$ has a submodule which is isomorphic to $H_i(W_{f_0};A)$ where they are seen as modules over $A$.
We can regard $e_{i,j}$ as an element of $H_i(M_0;A)$ by mapping this by the homomorphism between the homology groups induced by the section of the trivial bundle over $W_{f_0}-{\rm Int}\ N(\partial W_{f_0}) \subset W_{f_0}$ and we have a submodule $H_i(W_{f_0};A)_{H_i(M_0;A)}$ isomorphic to $H_i(W_{f_0};A)$. Consider the Poincar\'e duals to this and ${\rm PD}({e_{i,j}}^{\ast, {\mathcal{E}}_{X,A}}) \in H_{m-i}(M_0;A)$ and the cup product for this (suitably ordered) pair. We have a generator of $H_{m}(M_0;A)$. 

The set $\{{\rm PD}({e_{i,j}}^{\ast, {\mathcal{E}}_{X,A}}) \in H_{m-i}(M_0;A)\}_{j \in J}$ generates a submodule $H_{n-i}(W_{f_0},\partial W_{f_0};A)_{H_{m-i}(M_0;A)} \subset H_{m-i}(M_0;A)$.
The cup product for the ordered pair
$({e_{i_1,j_1}}^{\ast, {\mathcal{E}}_{X,A}},{e_{i_2,j_2}}^{\ast, {\mathcal{E}}_{X,A}})$ such that the sum of the degrees is $m$ is the zero element. This is due to the fact that the sum of the degrees of ${\rm PD}({e_{i_1,j_1}}^{\ast, {\mathcal{E}}_{X,A}})$ and ${\rm PD}({e_{i_2,j_2}}^{\ast, {\mathcal{E}}_{X,A}})$ is $$(i_1-(m-n))-(i_2-(m-n))=(i_1+i_2)-(m-n)-(m-n)=m-2(m-n)=2n-m<2n-n=n$$ in the case where they are seen as elements of $H_n(W_{f_0},\partial W_{f_0};A)$ and Poincar\'e duality theorem or intersection theory for $M_0$ and $W_{f_0}$. Such theorem and theory are also essential in the proof of (\ref{thm:3.3.3}).

This completes the proof of  (\ref{thm:3.3.1}), (\ref{thm:3.3.2}), (\ref{thm:3.3.3}), and (\ref{thm:3.3.4}).

Note again that Poincar\'e duality or intersection theory presents key tools here. 
This is also discussed in the proof of Theorem \ref{thm:1} in \cite{kitazawa8}. We do not know whether we can have more precise and concrete arguments on Poincar\'e duals as in \cite{kitazawa8} since the situation is much more generalized than the original one. However, this is conjectured by us that we can do in considerable scenes here.

We show remaining properties (\ref{thm:3.3.5}) and (\ref{thm:3.3.6}), which are new. 

We can see that as a submodule of the homology group $H_{n-1}(W_{f_0};A)$, we have some free module of rank $l_2-1$ over $A$ and its basis each of whose element is represented by the ($n-1$)-dimensional closed and connected manifold $(W_{f_0}-{\rm Int}\ N(\partial W_{f_0})) \bigcap \partial N(C_j)$ with $1 \leq j \leq l_2-1$. 
We explain about this via Poincar\'e duality and intersection theory for $W_{f_0}$. Poincar\'e duality theorem implies that $H_{n-1}(W_{f_0};A)$ is isomorphic to $H^1(W_{f_0}, \partial W_{f_0};A)$. We can choose a basis $\{e_{0,1,j}\}_{j=1}^{l_2-1}$ of a suitable submodule of $H_1(W_{f_0}, \partial W_{f_0};A)$ and each of the elements of the basis is represented by a smoothly embedded closed interval $I_j$ in $W_{f_0}$ satisfying the following conditions.
\begin{itemize}
	\item Each element $e_{0,1,j}$ is represented by a smoothly embedded closed interval $I_j$ whose interior is embedded in the interior of $W_{f_0}$ and each of whose two points in the boundary $\partial I_j$ is embedded in the connected component $C_0$ and the connected component $C_j$, respectively. For each $e_{0,1,j}$, the image $I_j$ of the embedding and $C_0$ intersect in a one-point set satisfying the transversality: the intersection of the image of the differential and the tangent bundle at the point of $C_0$ is the trivial vector space. For each $e_{0,1,j}$, the image $I_j$ of the embedding and $C_j$ intersect in a one-point set satisfying the transversality: the intersection of the image of the differential and the tangent bundle at the point of $C_j$ is the trivial vector space.
	
	\item For each $e_{0,1,j}$, the image $I_j$ of the embedding and $(W_{f_0}-{\rm Int}\ \partial W_{f_0}) \bigcap \partial N(C_j)$ intersect in a one-point set satisfying the transversality: the intersection of the images of the differentials at these two points is the trivial vector space. The image $I_j$ of the embedding and $(W_{f_0}-{\rm Int}\ \partial W_{f_0}) \bigcap \partial N(C_0)$ intersect in one-point set satisfying the similar transversality.
	\item The $l_2-1$ smoothly embedded intervals in $\{I_j\}$ are mutually disjoint.

\end{itemize}

This also explains about the following two.
\begin{itemize}
\item The desired fact that as a submodule of the homology group $H_{n-1}(W_{f_0};A)$, we have some free submodule of rank $l_2-1$ over $A$ and its basis each of whose element is represented by the ($n-1$)-dimensional closed and connected manifold $(W_{f_0}-{\rm Int}\ N(\partial W_{f_0})) \bigcap \partial N(C_j)$ with $1 \leq j \leq l_2-1$. LEMMA 3.1 of \cite{nishioka} also says that
in the case $H_1(W_{f_0};A)$ is the trivial group, the submodule can be chosen as the original module $H_{n-1}(W_{f_0};A)$.
\item Furthermore, $H_{n-1}(\partial W_{f_0};A)$ is mapped into $H_{n-1}(W_{f_0};A)$ as a monomorphism.
\end{itemize}

In addition, by considering the local structures of the maps, we may also say that the restrictions of the map $q_f$ to the preimages of the intervals give the compositions of the projections of trivial smooth bundles whose fibers are products of standard spheres with height functions on the unit disks. We study related functions later.

We consider a Mayer-Vietoris sequence
$$\rightarrow H_i(M_0;A) \rightarrow H_{i-1}(P_{f_0} \bigcap B_{C,f_0};A) \rightarrow $$
$$H_{i-1}(P_{f_0};A) \oplus H_{i-1}(B_{C,f_0};A) \rightarrow H_{i-1}(M_0;A) \rightarrow $$
and let $i>n$ here.
We consider a subbundle $S_{P,B,C_j,f_0}$ of the trivial smooth bundle defined by the map $q_{f_0}:M_0 \rightarrow W_{f_0}$ over $(W_{f_0}-{\rm Int}\ N(\partial W_{f_0})) \bigcap \partial N(C_j)$ whose fiber is the product of $l_1$ manifolds satisfying the following rules. They are also in our related exposition in our present theorem for each $1 \leq j \leq l_2-1$.

\begin{itemize}
	\item The fiber is represented as the product ${\prod}_{j^{\prime}=1}^{l_1} S_{P,B,C_j,f_0,j^{\prime}}$.
	\item $S_{P,B,C_j,f_0,j^{\prime}}:=S^{m_{l_1}(j^{\prime})}$ if $j^{\prime} \notin m_{{\mathcal{C}}_{\bar{N}},l_1}(\{C_0,C_j\})$.
	\item $S_{P,B,C_j,f_0,j^{\prime}}$ is a one-point set in $S^{m_{l_1}(j^{\prime})}$ if $j^{\prime} \in m_{{\mathcal{C}}_{\bar{N}},l_1}(\{C_0,C_j\})$.
\end{itemize}

The total space of this bundle is also a smooth compact and connected submanifold with no boundary in a closed and connected manifold $P_{f_0} \bigcap B_{C_j,f_0}$. Consider each element of $H_{i-1}(P_{f_0} \bigcap B_{C,f_0};A)$ represented by the total space of a subbundle ${S^{\prime}}_{P,B,C_j,f_0,a}$ of $S_{P,B,C_j,f_0}$ over $(W_{f_0}-{\rm Int}\ N(\partial W_{f_0})) \bigcap \partial N(C_j)$ whose fiber is represented by the product of $l_1$ manifolds satisfying the following rules as in our related exposition in our present theorem. 
\begin{itemize}
	\item The fiber is represented as the product ${\prod}_{j^{\prime}=1}^{l_1} {S^{\prime}}_{P,B,C_j,f_0,j^{\prime},a}$.
\item ${S^{\prime}}_{P,B,C_j,f_0,j^{\prime},a}$ is ${S}_{P,B,C_j,f_0,j^{\prime}}$ or a one-point set there.
\item The fiber is not a one-point set. It is ($i-n$)-dimensional.

\end{itemize}

This element, denoted by $[{S^{\prime}}_{P,B,C_j,f_0,a}] \in H_{i-1}(P_{f_0} \bigcap B_{\mathcal{C},f_0};A)$, is mapped both into the two direct summands in the Mayer-Vietoris sequence and the resulting values are not the zero elements. Let the resulting element denoted by $([{S^{\prime}}_{P,B,C_j,f_0,a}] \in H_{i-1}(P_{f_0};A),[{S^{\prime}}_{P,B,C_j,f_0,a}] \in H_{i-1}(B_{\mathcal{C},f_0};A))$. We consider the element of the form $([{S^{\prime}}_{P,B,C_j,f_0,a}],0_A) \in H_{i-1}(P_{f_0};A) \oplus H_{i-1}(B_{C,f_0};A)$ where $0_A$ is for the zero element. This is not in the image of the homomorphism in the Mayer-Vietoris sequence. We explain about this.

We consider an arbitrary element of the form $([{S^{\prime}}_{P,B,C_j,f_0,a}],c) \in H_{i-1}(P_{f_0};A) \oplus H_{i-1}(B_{C,f_0};A)$. 
Just before, we have said that $H_{n-1}(\partial W_{f_0};A)$ is mapped into $H_{n-1}(W_{f_0};A)$ as a monomorphism. We can see that
$([{S^{\prime}}_{P,B,C_j,f_0,a}],c) \in H_{i-1}(P_{f_0};A) \oplus H_{i-1}(B_{C,f_0};A)$ must be mapped from $[{S^{\prime}}_{P,B,C_j,f_0,a}] \in H_{i-1}(P_{f_0} \bigcap B_{\mathcal{C},f_0};A)$. In addition, $c$ is, by the definition, shown to be not $0_A$.
We present arguments to understand this.

This is understood by Künneth formula, fundamental tools in knowing the (co)homology groups and the cohomology rings of the products of finitely many spaces, for example. Remember that $P_{f_0}$, $B_{C_j,f_0}$ and $P_{f_0} \bigcap B_{C_j,f_0}$ are the total spaces of naturally obtained trivial bundles. For this, note also that the fiber of the trivial bundle $P_{f_0} \bigcap B_{C_j,f_0} \subset B_{C_j,f_0}$ is the boundary of the fiber of $B_{C_j,f_0}$ obtained by replacing the exactly one unit disk in the product for the fiber by the unit sphere being the boundary of the previous unique unit disk and that the former bundle is a subbundle of the latter bundle. Remember also the conditions on the unit disks.

These arguments show that $([{S^{\prime}}_{P,B,C_j,f_0,a}],0_A) \in H_{i-1}(P_{f_0};A) \oplus H_{i-1}(B_{C,f_0};A)$ is not in the image of the homomorphism in the Mayer-Vietoris sequence. We have a submodule of $H_{i-1}(M_0;A)$ generated by
 all elements obtained by mapping the elements $([{S^{\prime}}_{P,B,C_j,f_0,a}],0_A) \in H_{i-1}(P_{f_0};A) \oplus H_{i-1}(B_{C,f_0};A)$ in the Mayer-Vietoris sequence. The element is also understood as the value of the homomorphism induced by the canonical inclusion at $[{S^{\prime}}_{P,B,C_j,f_0,a}] \in H_{i-1}(P_{f_0};A)$.

We investigate the restriction of the map $q_{f_0}$ to ${q_{f_0}}^{-1}(I_j)$ and define this map as a map onto $I_j$. We compose a smooth embedding into $\mathbb{R}$. By suitably restricting this function, we have the composition of the projection of a trivial smooth bundle whose fiber is diffeomorphic to ${\prod}_{j^{\prime}=1}^{l_1} S_{P,B,C_j,f_0,j^{\prime}}$ over a unit sphere with a Morse-Bott function in the following. Let this function obtained by this restriction denoted by $q_{f_0,I_j}$.
\begin{itemize}
	\item If the set $m_{{\mathcal{C}}_{\bar{N}},l_1}(\{C_0,C_j\})$ consists of exactly one element $j_0$, then the unit sphere of the base space is $S^{m_{l_1}(j_0)+1}$ and the Morse-Bott function is Morse and a height function of the unit sphere.
	\item If the set $m_{{\mathcal{C}}_{\bar{N}},l_1}(\{C_0,C_j\})$ consists of exactly two elements $j_{0,1}$ and $j_{0,2}$, then the unit sphere is $S^{m_{l_1}(j_{0,1})+m_{l_1}(j_{0,2})+1}$.
\end{itemize}	
	We explain about remaining results to show. We apply Poincar\'e duality and intersection theory. 
In short, ${S^{\prime}}_{P,B,C_j,f_0,a}$ and the manifold of the domain of a suitable restriction of the function $q_{f_0,I_j}$ intersect in a one-point set satisfying some transversality as presented. Moreover, we need to choose the restriction so that the projection onto the unit sphere is the projection of the subbundle whose fiber is defined as follows.

\begin{itemize}
	\item The fiber is represented as the product ${\prod}_{j^{\prime}=1}^{l_1} {S^{\prime \prime}}_{P,B,C_j,f_0,j^{\prime},a}$.
	\item ${S^{\prime \prime}}_{P,B,C_j,f_0,j^{\prime},a}$ is ${S}_{P,B,C_j,f_0,j^{\prime}}$ 
	if both the following two are satisfied.
	\begin{itemize}
		\item ${S^{\prime }}_{P,B,C_j,f_0,j^{\prime},a}$ is defined to be a one-point set.
		\item ${S}_{P,B,C_j,f_0,j^{\prime}}$ is not defined to be a one-point set and it is defined to be a standard sphere.
		\end{itemize}
	\item ${S^{\prime \prime}}_{P,B,C_j,f_0,j^{\prime},a}$ is a one-point set if either of the following two is satisfied.
	\begin{itemize}
		\item ${S^{\prime }}_{P,B,C_j,f_0,j^{\prime},a}$ is defined to be ${S}_{P,B,C_j,f_0,j^{\prime}}$ and this is a standard sphere.
		\item ${S}_{P,B,C_j,f_0,j^{\prime}}$ is defined to be a one-point set.
		\end{itemize}
\end{itemize}

The images in the family $\{I_j\}$ are assumed to be mutually disjoint. These arguments on intersection theory with the previous arguments on the submodule of $H_{n-1}(W_{f_0};A)$ and the Mayer-Vietoris sequence with K\"unneth theorem for the total spaces of the trivial bundles prove (\ref{thm:3.3.5}) for example.

For (\ref{thm:3.3.6}), the preimage of $W_{f_0}-{\rm Int}\ N(\partial W_{f_0})$ by $q_{f_0}$, which is (the total space of) a trivial smooth bundle, admits a section and it is important. This section is also important in defining the submodule of (\ref{thm:3.3.2}), as presented in the related argument, for example. More precisely, by the structures of our maps, consisting of trivial smooth bundles, we can see that ${S^{\prime}}_{P,B,C_j,f_0,a}$ and the image of the section can be moved by isotopies in the (piecewise) smooth category to be mutually disjoint. This with related explicit intersection theory is a main ingredient to show (\ref{thm:3.3.6}).

This completes the proof. 
\end{proof}
\begin{Ex}
	\label{ex:2}
\begin{enumerate}
\item \label{ex:2.1} Theorem \ref{thm:1} is a specific case of Theorem \ref{thm:3} of course. Submodules in (\ref{thm:3.3.5}) are trivial. On Theorem \ref{thm:1}, we do not consider generalizations of the property (\ref{thm:1.1}). This seems to be difficult.
\item \label{ex:2.2}
In the situation of Theorem \ref{thm:3}, let $\bar{N}:=S^{n-1} \times D^1$ and let it be smoothly embedded in $N:={\mathbb{R}}^n$. 
Here
$l_2=2$ and let the two connected components of the boundary denoted by $C_0$ and $C_1$, respectively. Suppose also that
$l_1=2$ and that $m_{l_1}(1)=m_2(1)=k_1>0$ and $m_2(2)=k_2>0$ hold. Suppose also that $m_{{\mathcal{C}}_{\bar{N}}}(C_j)=j+1$ for $j=0,1$.

Then we have a simply generalized special generic map $f_0:M_0:=S^{k_1+k_2+1} \times S^{n-1} \rightarrow {\mathbb{R}}^n$ obtained in the following way. 
\begin{itemize}
	\item First consider the product map of a Morse-Bott function on $S^{k_1+k_2+1}$ with exactly two singular values and the identity map on $S^{n-1}$. 
For the Morse-Bott function, let the preimage of each point in the image diffeomorphic to $S^{k_1} \times S^{k_2}$.
	\item Embed the manifold of the target before into $N={\mathbb{R}}^n$.
\end{itemize}  
This is also for Theorem \ref{thm:3} where submodules in (\ref{thm:3.3.5}) are trivial. 
\item \label{ex:2.3}
In the situation of Theorem \ref{thm:3}, suppose that
the map $m_{{\mathcal{C}}_{\bar{N}}}$ is constant and that $l_1 \geq 1$. In the case $l_1=1$, this is for Theorem \ref{thm:1}. In the case $l_1>1$, we can have our desired map $f_0$ as the composition of a projection of a trivial smooth bundle over a closed and connected manifold with a special generic map there as in Theorem \ref{thm:1} or Theorem \ref{thm:3}.

By this explicit theory or more general one, we can construct simply generalized special generic maps on $m$-dimensional manifolds into ${\mathbb{R}}^n$ where the manifolds do not admit special generic maps into ${\mathbb{R}}^n$. 

For example, we have a simply generalized special generic map on $S^2 \times S^2 \times S^2$ into ${\mathbb{R}}^3$ as the composition of the projection of a trivial smooth bundle over $S^2 \times S^2$ whose fiber is $S^2$ with a special generic map into ${\mathbb{R}^3}$ of Theorem \ref{thm:2} (\ref{thm:2.2}) and (\ref{thm:2.3}) or Example \ref{ex:1}. The non-existence of special generic maps in this case follows from Theorem \ref{thm:2} (\ref{thm:2.2}) and we can also know from the preprint \cite{kitazawa1} of the author. 

We investigate this in terminologies of Theorem \ref{thm:3} more precisely. We can see this is for $(m,l_1,l_2)=(6,2,2)$.
We can take $(i_1,i_2)=(2,5)$ and $i=i_2=5$ in Theorem \ref{thm:3} and the rank of the free submodule of the module $H_{i_2-1}(M_0;A)=H_{4}(M_0;A)$ in (\ref{thm:3.3.5}) is $1$. We do not cover the free submodule of rank $1$ generated by an element represented by a submanifold regarded as the total space $S^2 \times S^2$ of the trivial bundle over $S^2$, embedded canonically in the manifold $S^2 \times S^2 \times S^2$. In the situation of Theorem \ref{thm:3}, for $i=5$, we can consider ${S^{\prime \prime}}_{P,B,C_1,f_0,1,a} \times {S^{\prime \prime}}_{P,B,C_1,f_0,2,a}$ and it can be chosen as a fiber of the trivial bundle $S^2 \times S^2 \times S^2$ over the total space of the trivial bundle $S^2 \times S^2$. We may regard $S_{P,B,C_1,f_0,1} \times S_{P,B,C_1,f_0,2} ={S^{\prime}}_{P,B,C_1,f_0,1,a} \times {S^{\prime}}_{P,B,C_1,f_0,2,a}$ in the theorem and these spaces can be also regarded as a fiber of the trivial bundle $S^2 \times S^2 \times S^2$ over the total space of the trivial bundle $S^2 \times S^2$. Last note that the set of all labels "$a$" is regarded as a one-element set.

\end{enumerate}
\end{Ex}

Investigating our method of our proof of Theorem \ref{thm:3} and generalizing Example \ref{ex:2} (\ref{ex:2.3}) for example, we also have another new result.

\begin{MainThm}
	\label{mthm:2}
	In the situation of Theorem \ref{thm:3}, suppose that the given map $m_{{\mathcal{C}}_{\bar{N}},l_1}$ is not surjective. Then we can have a desired map $f_0:M_0 \rightarrow N$ as a map represented as the composition of the projection of a trivial smooth bundle whose fiber is the product of finitely many standard spheres over another closed and connected manifold $M_{0,0}$ with any simply generalized special generic map $f_{0,0}:M_{0,0} \rightarrow N$ enjoying the following two.
\begin{enumerate}
\item $f_{0,0}:M_{0,0} \rightarrow N$ is obtained by Theorem \ref{thm:3}.
\item Furthermore, $f_{0,0}:M_{0,0} \rightarrow N$ is obtained by considering another suitable surjective map $m_{{\mathcal{C}}_{\bar{N}},{l_1}^{\prime}}$ with a suitable positive integer ${l_1}^{\prime}<l_1$ and another suitable map $m_{{l_1}^{\prime}}$ where notation is suitably and naturally considered and abused.
\end{enumerate}
\end{MainThm}

	\end{document}